\documentclass[11pt,twosided]{article}
\usepackage{dsfont}
\usepackage{bm}
\usepackage{float}
\usepackage{srcltx}
\usepackage{mathrsfs}
\usepackage{graphicx}
\usepackage{subfigure}
\usepackage{amsfonts}
\usepackage{indentfirst}
\usepackage{amssymb,amsfonts,amsmath, psfrag,eepic,colordvi}
\usepackage{diagbox}
\usepackage{caption}

\usepackage{leftidx}
\usepackage{pgf,tikz}
\usepackage{cite}
\usetikzlibrary{shapes,arrows}
\allowdisplaybreaks[4]
\usepackage{color}

\usepackage{amssymb,color}

\definecolor{c20}{rgb}{0.,0.7,0.}
\definecolor{c30}{rgb}{0.,0.,1.}
\definecolor{c40}{rgb}{1,0.1,0.7}
\definecolor{c50}{rgb}{1,0,0}
\definecolor{c60}{rgb}{0,0.9,0.1}

\usepackage{amssymb,color}

\newcommand{\BQN}{\begin{eqnarray}}
\newcommand{\EQN}{\end{eqnarray}}
\newcommand{\BQNY}{\begin{eqnarray*}}
\newcommand{\EQNY}{\end{eqnarray*}}

\newcommand{\BS}{\begin{sat}}
\newcommand{\ES}{\end{sat}}
\newcommand{\BT}{\begin{theo}}
\newcommand{\ET}{\end{theo}}
\newcommand{\BK}{\begin{korr}}
\newcommand{\EK}{\end{korr}}

\newcommand{\BD}{\begin{de}}
\newcommand{\ED}{\end{de}}

\newcommand{\BIT}{\begin{itemize}}
\newcommand{\EIT}{\end{itemize}}
\newcommand{\BDI}{\begin{description}}
\newcommand{\EDI}{\end{description}}

\newcommand{\BRM}{\begin{remarks}}
\newcommand{\ERM}{\end{remarks}}

\newcommand{\BEL}{\begin{lem}}
\newcommand{\EEL}{\end{lem}}

\newtheorem{theo}{Theorem}[section]

\newtheorem{lem}[theo]{Lemma}

\newtheorem{re}[theo]{Remark}
\newtheorem{remarks}[theo]{Remarks}

\newcommand{\COM}[1]{}

\begin{document}
\title{Number of claims and ruin time for a refracted risk process}
\author{Yanhong Li\thanks{Sichuan University, Chengdu, 610065, China.
Email: yanhonglink@qq.com}
\ \ Zbigniew Palmowski
\thanks{Faculty of Pure and Applied Mathematics,
Wroc\l aw University of Science and Technology,
Wyb. Wyspia\'nskiego 27, 50-370 Wroc\l aw, Poland, e-mail: zbigniew.palmowski@gmail.com}
\ \  Chunming Zhao\thanks{Department of Statistics,
School of Mathematics, Southwest Jiaotong University, Chengdu, Sichuan, 611756, China. Email:cmzhao@swjtu.cn}
 \ \ Chunsheng Zhang\thanks{School of mathematical Sciences and LPMC Nankai University, Tianjin 300071, China.
Email: zhangcs@nankai.edu.cn}
 }
\maketitle
\begin{abstract}
In this paper, we consider a classical risk model refracted at given level.
We give an  explicit expression for the joint density of the ruin time and the cumulative number of claims counted up to ruin time.
The proof is based on solving some integro-differential equations and employing the Lagrange's Expansion Theorem.
\end{abstract}


\section{Introduction}
The joint density of the ruin time and the numbers of claims counted until ruin time has been already studied for a classical risk process over last years.
Dickson \cite{Di2012} derived special expression for it using probabilistic arguments. Landriault et al. \cite{LaShWi2011} analyzed this object
for the Sparre Andersen risk model with the exponential claims. Later Frostig et al. \cite{FrPiPo2012} generalized it to the case of a renewal risk model
with the phase-type claims and inter-arrival times. The main tool used there was the duality between the risk model and a workload of a single server queueing
model. Zhao and Zhang \cite{ZhZh2013} considered a delayed renewal risk model, where the claim size is Erlang$(n)$ distributed and the inter-arrival time is assumed to be infinitely divisible.

Our goal is to derive expression for the joint density of the ruin time and the numbers of claims counted until ruin time for a refracted classical risk process (see
Kyrianou and Loeffen \cite{KyLo2010} for a formal definition). It is also called a compound Poisson risk model under a threshold strategy. The latter process is a classical risk process whose dynamic is changed by subtracting off a fixed linear drift
whenever the cumulative risk process is above a pre-specified level $b$. This subtracting of the linear drift corresponds to the dividend payments
and the considered strategy is also known as a threshold strategy.
Dividend strategies for insurance risk models were first proposed by De Finetti \cite{De1957} to reflect more realistically the surplus cash flows in an insurance portfolio.
More recently, many kind of risk related quantities under threshold dividend
strategies have been studied by Lin and Pavlova \cite{LiPa2006}, Zhu and Yang \cite{ZhYa2008},
Lu and Li \cite{LiLu2009}, \cite{LiLu2016}, \cite{LiLu2},
Badescu, Drekic and Landriault \cite{BaDrLa2007}, Gao and Yin \cite{GaYi2008} (see references therein).
The case when the drift of the refracted process is disappearing (everything above threshold $b$ is paid as dividends) is called barrier strategy, see
Lin et al. \cite{LiWiDr2003}, Li and Garrido \cite{LiGa2004a}, Zhou \cite{Zh2005} and in the references therein.

The paper is organized as follows. In Section 2 we define the model we deal with in this paper.
In Section 3 we recall properties of the translation operator and the root of the Lundberg fundamental equation. In particular, we introduce the Lagrange's  expansion theorem and some notation.  In Section 4 we construct two integro-differential equations identifying the joint Laplace transform of joint density of the numbers of claims counted up to ruin time
and the ruin time.
Analytical solutions of these two integro-differential equations are given in Section 5. Applying the Lagrange's  expansion theorem in Section 6
we give the expression for above mentioned density.

\section{Model}
The classical risk process is given by
\begin{equation}\label{CL}
U(t)=u+c_1t-S(t),
\end{equation}
where $U(0)=u$ denotes initial capital, $c$ is the premium rate and
$S(t)=\sum_{i=1}^{N_t}X_i$ represents the total amount of claims appeared up to time $t\geq 0$. That is,  $\{X_i\}_{\{i\in \mathbb{N}\}}$ are non-negative i.i.d. random variables with pdf $f(x)$ and cdf $F(x)$ and $\{N_t\}_{\{t\geq 0\}}$ is an independent Poisson process with a parameter $\lambda$.
To take into account dividend payments paid when regulated process (after deduction of dividends) is above fixed threshold level $b>0$, we consider so-called refracted process given formally for $c_2<c_1$ by:
\begin{equation}
dU_b(t)=\left\{
  \begin{array}{lll}
&c_1dt-dS(t),&0\leq U_b(t)\leq b\\
&c_2dt-dS(t),& U_b(t)>b
  \end{array}
\right.
\end{equation}
and $U_b(0)=u$.
In this case $c_1-c_2$ denotes intensity of dividend payments, see Figure 1.

\begin{center}
\definecolor{qqttff}{rgb}{0,0.2,1}
\begin{tikzpicture}[line cap=round,line join=round,>=triangle 45,x=0.5cm,y=0.5cm]
\clip(-2.65,-5) rectangle (22.64,14.03);
\draw (0.73,12.51) node[anchor=north west] {$U_b(t)$};
\draw (0,4)-- (2,6);
\draw [dotted] (2,6)-- (2,5);
\draw (2,5)-- (3,6);
\draw [dotted] (3,6)-- (3,3);
\draw (3,3)-- (7,7);
\draw (7,7)-- (9,8);
\draw [dotted] (9,8)-- (9,4);
\draw (9,4)-- (12,7);
\draw (12,7)-- (13,7.5);
\draw [dotted] (13,7.5)-- (13,2);
\draw (13,2)-- (14,3);
\draw [dotted] (14,3)-- (14,-2);
\draw [->] (0,-2.79) -- (0,12);
\draw [->] (-2,0) -- (21,0);
\draw [color=qqttff](0,7)-- (21,7);
\draw (-0.74,8.00) node[anchor=north west] {$b$};
\draw (-0.87,4.74) node[anchor=north west] {$u$};
\draw (-0.87,0.04) node[anchor=north west] {$0$};
\draw (20.55,0.04) node[anchor=north west] {$t$};
\draw (13.30,0.04) node[anchor=north west] {$\tau$};
\draw (5,5) node[anchor=north west] {$c_1$};
\draw (7.5,8.7) node[anchor=north west] {$c_2$};
\draw (-0.77,-3.51) node[anchor=north west] {${\rm Figure\  1:  \; Graphical   \; representation  \; of \;  the \;  surplus  \; process}  \; U_b(t).$};
\end{tikzpicture}
\end{center}

Throughout this paper, we will assume  that $c_2>\lambda EX_1$, which means refracted process $U_b(t)$ tends to infinity almost surely.
We can then consider the ruin time:
\begin{equation*}
    \tau=\inf\{t>0,U_b(t)<0\},
\end{equation*}
($\tau􏰽􏱨=\infty$ if ruin does not occur). Note that $N_\tau$ represents the number of claims counted until the ruin time.
The main goal of this paper is identification of the density of $(\tau, N_\tau)$.
We start from analyzing its Laplace transform:
\begin{eqnarray}\label{phi}
    \phi(u)&=&E[r^{N_\tau}e^{-\delta\tau}\mathbb{I}(\tau<\infty)|U_b(0)=u]\\
           &=&\sum_{n=1}^\infty r^n\int_0^\infty e^{-\delta t}w(u,n,t)dt,\label{define_phi}
\end{eqnarray}
where
$$w(u,n,t)=P(N_\tau=n, \tau\in dt|U_b(0)=u)/dt$$
is the joint density of $(\tau,N_\tau)$ when $U_b(0)=u$. In above definition we have $\delta>0$ and $r\in(0,1]$.
Later we will use the following notation
\begin{equation}\label{w1def}
w_1(u,n,t)=w(u,n,t) \qquad\text{for $u\leq b$}
\end{equation}
and
\begin{equation}\label{w1def}
w_2(u,n,t)=w(u,n,t) \qquad\text{for $u> b$.}
\end{equation}

\section{Preliminaries}
In this section we introduce few facts used further in this paper.
We start from recalling the translation operator $T_s$; see Dickson and Hipp \cite{DiHi2001}.
For any integrable real-valued function $f$ it is defined as
\begin{equation*}
T_s f(x)=\int_x^{\infty}{e^{-s(y-x)}f(y)}dy,\ \   x\geq0.
\end{equation*}
The operator $T_{s}$ satisfies the following properties:
\begin{enumerate}
  \item $ T_s f(0)=\int_0^{\infty}{e^{-s x}f(x)dx}=\hat f(s)$ which is the Laplcae transform of $f$;
\item The operator $T_s$ is commutative, i.e. $T_sT_r=T_rT_s$. Moreover, for $s\neq r$ and $x\geq0$
\begin{equation}\label{Dickson's Operator}
T_{s}T_{r}f(x)=T_{r}T_{s}f(x)=\frac{T_s f(x)-T_r f(x)}{r-s}.
\end{equation}
\end{enumerate}
More properties of the translation operator $T_s$ can be found in Li and Garrido \cite{LiGa2004b} and Gerber and Shiu \cite{GeSh2005}.

For any function $g$ we will denote by $\hat{g}(s)$ its Laplace Transform, that is $\hat{g}(s)=\int_0^\infty e^{-s x}g(x)\;dx$.
Next, for $i=1,2$  let $\rho_{i}$ be the positive root of the Lundberg fundamental equation
\begin{equation}\label{Lundberg_fundamental_equation}
    c_{i}s-(\lambda+\delta)+\lambda r\hat{f}(s)=0.
\end{equation}
The positive roots always exists for $\delta>0$; see Figure 2.\\
\begin{figure}
\centering
\includegraphics[width=10cm]{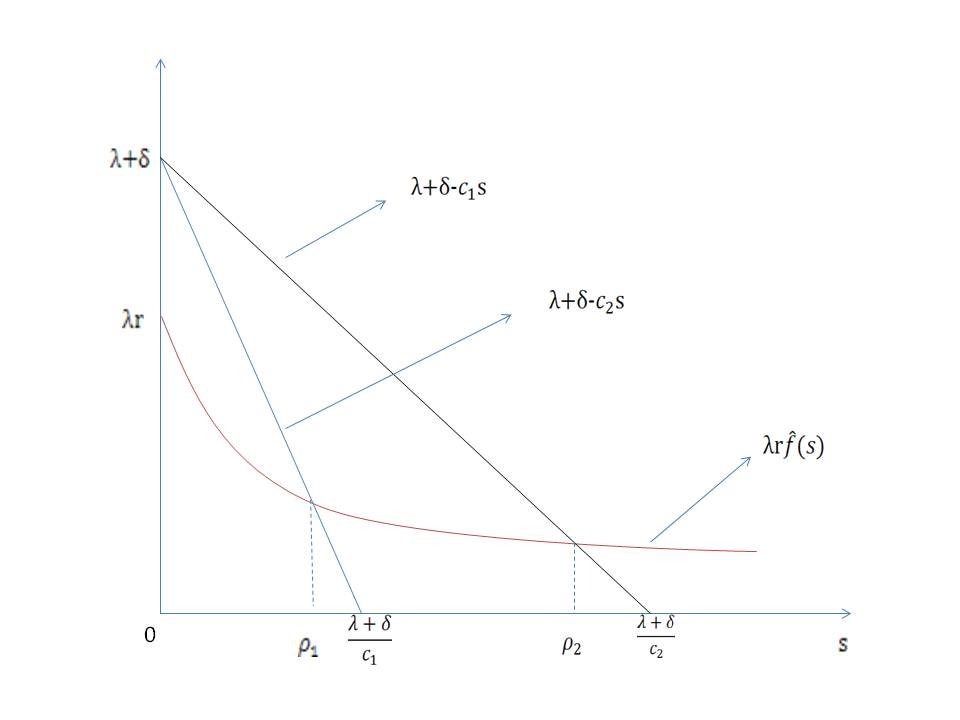}\label{figure_LFE}
\caption*{Figure 2: Roots for Lundberg's fundamental equation.}
\label{figure_LFE}
\end{figure}

\textbf{Lagrange's Expansion Theorem.}\label{LIFT}
In this paper we will also use the Lagrange's Expansion Theorem; see pages 251-326 of Lagrange \cite{La1770}.
Given two functions $\alpha(z)$ and $\beta(z)$ which are both analytic on and inside a contour $D$ surrounding a point $a$, if $r$ satisfies the inequality
\begin{equation}
|r \beta(z)|<|z-a|,
\end{equation}
for every $z$ on the perimeter of $D$, then $z-a-r \varphi(z)$, as a function of $z$, has exactly one zero $\eta$ in the interior of $D$,  and we further have:
\begin{equation}
    \alpha(\eta)=\alpha(a)+\sum_{k=1}^{\infty}\frac{r^k}{k!}\frac{d^{k-1}}{dx^{k-1}}\big(\alpha'(x)\beta^k(x)\big)\!\big|_{x=a}.
\end{equation}

Finally, we define also the impulse function
\begin{equation*}
\delta_x(t)=\left\{
  \begin{array}{lll}
0,&&t\neq x\\
\infty,&&t=x
  \end{array}
\right.
\end{equation*}
with $\int_0^\infty\delta_x(t)dt=1$. We denote $g^{k*}, k\ge 0$, with $g^{1*}=g$ and $g^{0*}(t)=\delta_0(t)$ the $k$-fold convolution of $g$ with itself, where
\begin{equation*}
(g*h)(t)=\int_{0}^{t}{g(x)h(t-x)dx},\  \ \ t\ge0
\end{equation*}
for two functions $g$ and $h$ supported on $(0,\infty)$.

\section{Integro-differential equations for the joint Laplace transform}
In this section, we derive two integro-differential equations identifying $\phi(u)$ defined in
(\ref{phi}). We will follow the idea given in Lin and Pavlova \cite{LiPa2006}.
Denote
\begin{equation}
\phi(u)=\left\{
  \begin{array}{lll}
\phi_1(u),&&u\leq b,\\
\phi_2(u),&&u>b.
  \end{array}
\right.
\end{equation}

\begin{theo}
The joint Laplace transform $\phi$ satisfies the following integro-differential equations:

\begin{equation}\label{phi'}
\left\{
  \begin{array}{lll}
\!\!\phi_1'(u)\!\!\!&=&\!\!\!\frac{\lambda+\delta}{c_1}\phi_1(u)-\frac{\lambda r}{c_1}\int_0^u\phi_1(u-x)f(x)dx-\frac{\lambda r}{c_1}\bar{F}(u),\ \ \ \ \ \ \ \ \ \ \ \ \ \ \ \ \ \ \ \ \ \ \ \ \ \ \ \ \ \ \ 0\leq u\leq b\\
\!\!\phi_2'(u)\!\!\!&=&\!\!\!\frac{\lambda +\delta}{c_2}\phi_2(u)-\frac{\lambda r}{c_2}\left(\int_0^{u-b}\phi_2(u-x)f(x)dx+\int_{u-b}^u\phi_1(u-x)f(x)dx\right)-\frac{\lambda r}{c_2}\bar{F}(u),\ u>b
  \end{array}
\right.
\end{equation}
with the boundary condition
\begin{equation}\label{bc}
\phi_1(b)=\phi_2(b):=\lim_{u\rightarrow b^{+}}\phi_2(u).\end{equation}
\end{theo}
\begin{re}
\rm Note that from the integro-differential equations (\ref{phi'}) follows that the joint Laplace transform with initial surplus above the barrier depends on the respective function with initial surplus below the barrier, but the reverse relationship does not hold true.
\end{re}

\textbf{Proof.}
Let first $0\leq u\leq b$. Then conditioning on the occurrence of the first claim we will have two cases: the first claim occurs before the surplus has reached the barrier level $b$
or it occurs after reaching this barrier. There are also two other cases at the moment of the arrival of the first claim:
either the risk process starts all over again with new initial surplus or the first claim leads already to ruin. Hence:
\begin{eqnarray}\label{phi_1_0}
\phi(u) &=&\phi_1(u)\nonumber\\
        &=&\int_0^{\frac{b-u}{c_1}}\lambda re^{-\lambda t}e^{-\delta t}\left(\int_0^{u+c_1t}\phi(u+c_1t-x)f(x)dx+\bar{F}(u+c_1t)\right)dt\nonumber\\
        &&+\int_{\frac{b-u}{c_1}}^\infty\lambda re^{-\lambda t}e^{-\delta t}\left(\int_0^{b+c_2(t-\frac{b-u}{c_1})}\phi(b+c_2(t-\frac{b-u}{c_1})-x)f(x)dx+\bar{F}(b+c_2(t-\frac{b-u}{c_1}))\right)dt\nonumber\\
        &=&\lambda r\int_0^{\frac{b-u}{c_1}}e^{-(\lambda+\delta)t}\gamma(u+c_1t)dt+\lambda r\int_{\frac{b-u}{c_1}}^\infty e^{-(\lambda+\delta)t}\gamma(b+c_2(t-\frac{b-u}{c_1})dt,
\end{eqnarray}
where $\gamma(t)=\int_0^t\phi(t-x)f(x)dx+\bar{F}(t)$. \\
Changing variables in (\ref{phi_1_0}) and rearranging leads to the following equation for $0\leq u\leq b$:
\begin{equation}\label{phi_1_1}
    \phi_1(u)=\frac{\lambda r}{c_1}e^{(\lambda+\delta)u/c_1}\int_u^be^{-(\lambda+\delta)t/c_1}\gamma(t)dt+\frac{\lambda r}{c_2}e^{(\lambda+\delta)u/c_1}\int_b^\infty e^{-(\lambda+\delta)[t-(c_1-c_2)b/c_1]/c_2}\gamma(t)dt.
\end{equation}
Differentiating both sides of (\ref{phi_1_1}) with respect to $u$ yields
first equation.

Similarly, for $u>b$ we have:
\begin{eqnarray}\label{phi_2_0}
\phi(u)&=&\phi_2(u)\nonumber\\
&=&\int_0^\infty\lambda re^{-\lambda t}e^{-\delta t}\left(\int_0^{u+c_2t}\phi(u+c_2t-x)f(x)dx+\bar{F}(u+c_2t)\right)dt\nonumber\\
&=&\lambda r\int_0^\infty e^{-(\lambda+\delta)t}\gamma(u+c_2t)dt\nonumber\\
&=&\frac{\lambda r}{c_2}e^{(\lambda+\delta)u/c_2}\int_u^\infty e^{-(\lambda+\delta)t/c_2}\gamma(t)dt.
\end{eqnarray}
Differentiating both sides of (\ref{phi_2_0}) with respect to $u$ produces the second equation.

Note also that from equations (\ref{phi_1_1}) and (\ref{phi_2_0}) it follows that $\phi(u)$ is continuous at $u=b$ and hence (\ref{bc}) holds.
This completes the proof.
\vspace{3mm} \hfill $\Box$

\section{The analytical expression for $\phi(u)$}
In this section, we derive the analytical expression for $\phi_i(u)$ ($i=1,2$) using the translation operator introduced in Section 3.

\begin{theo}
The function $\phi_2(u)$ can be expressed analytically as follows:
\begin{equation}\label{renewal equation 2}
   \phi_2(u)=\sum_{n=0}^{\infty}\left(\frac{\lambda r}{c_2}\right)^{n+1} (T_{\rho_2}f)^{n*}*h(u-b),\ \ \ \ u>b,
\end{equation}
where
\begin{equation}\label{h}
    h(u):=\int_{u}^{u+b}\phi_1(u+b-x)T_{\rho_2}f(x)dx+T_{\rho_2}\bar{F}(u+b).
\end{equation}
\end{theo}
\textbf{Proof.} We adopt the approach of Willmot and Dickson \cite{WiDi2003}.
Consider the second equation in (\ref{phi'}) for $u>b$.
For a fixed $s>0$, we multiply both sides of this
equation by $e^{-s(u-b)}$ and integrate it with respect to $u$ from $b$ to $\infty$:
\begin{eqnarray*}
   c_2\int_{b}^{\infty}e^{-s(u-b)}\phi_2'(u)du
         &=& (\lambda +\delta)T_s\phi_2(b)-\lambda r\int_{b}^{\infty}e^{-s(u-b)}\int_0^{u-b}\phi_2(u-x)f(x)dx du\\
         && -\lambda r\int_{b}^{\infty}e^{-s(u-b)}\int_{0}^b\phi_1(y)f(u-y)dy du-\lambda r T_s\bar{F}(b) \\
         &=& (\lambda +\delta)T_s\phi_2(b)-\lambda r\int_{0}^{\infty}e^{-sx}f(x)\int_{x+b}^{\infty}e^{-s(u-x-b)}\phi_2(u-x)du dx\\
         && -\lambda r\int_{0}^{b}\phi_1(y)\int_{b}^{\infty}e^{-s(u-b)}f(u-y)du dy-\lambda r T_s\bar{F}(b) \\
         &=& (\lambda +\delta)T_s\phi_2(b)-\lambda r\hat{f}(s)T_s\phi_2(b)-\lambda r\int_{0}^{b}\phi_1(y)T_sf(b-y)dy-\lambda r T_s\bar{F}(b).
\end{eqnarray*}
Integrating by parts gives:
\begin{equation*}
  c_2\int_{b}^{\infty}e^{-s(u-b)}\phi_2'(u)du = c_2sT_s\phi_2(b)-c_2\phi_2(b).
\end{equation*}
Hence
\begin{equation*}
   c_2sT_s\phi_2(b)-c_2\phi_2(b)=(\lambda +\delta)T_s\phi_2(b)-\lambda r\hat{f}(s)T_s\phi_2(b)-\lambda r\int_{0}^{b}\phi_1(y)T_sf(b-y)dy-\lambda r T_s\bar{F}(b)
\end{equation*}
and simple rearranging leads to:
\begin{equation}\label{phi 2}
   (c_2s-(\lambda+\delta)+\lambda r\hat{f}(s))T_s\phi_2(b)=c_2\phi_2(b)-\lambda r\int_{0}^{b}\phi_1(y)T_sf(b-y)dy-\lambda r T_s\bar{F}(b).
\end{equation}
Taking $s=\rho_2$ for the solution $\rho_2$ of the Lundberg Fundamental Equation (\ref{Lundberg_fundamental_equation}) gives
\begin{equation*}
    c_2\phi_2(b)=\lambda r\int_{0}^{b}\phi_1(y)T_{\rho_2}f(b-y)dy+\lambda r T_{\rho_2}\bar{F}(b).
\end{equation*}
Then equation (\ref{phi 2}) is equivalent to:
\begin{equation*}
    [c_2(s-\rho_2)+\lambda r\hat{f}(s)-\lambda r\hat{f}(\rho_2)]T_s\phi_2(b)=\lambda r\int_{0}^{b}\phi_1(y)[T_{\rho_2}f(b-y)-T_{s}f(b-y)]dy+\lambda r [T_{\rho_2}\bar{F}(b)-T_{s}\bar{F}(b)].
\end{equation*}
Now dividing above equation by $s-\rho_2$ and using property 2 of the translation operator introduced in Section 2 produces:
\begin{equation}\label{equation for phi2}
c_2T_s\phi_2(b)=\lambda rT_sT_{\rho_2}f(0)T_s\phi_2(b)+\lambda r\int_{0}^{b}\phi_1(y)T_sT_{\rho_2}f(b-y)dy+\lambda r T_sT_{\rho_2}\bar{F}(b).
\end{equation}
Inverting the translation operators of (\ref{equation for phi2}) yields the following renewal equation for $\phi_2(u)$:
\begin{equation}\label{phi 2 renewal equation}
    \phi_2(u)=\frac{\lambda r}{c_2}\left[\int_0^{u-b}\phi_2(u-x)T_{\rho_2}f(x)dx+\int_{u-b}^u\phi_1(u-x)T_{\rho_2}f(x)dx+T_{\rho_2}\bar{F}(u)\right].
\end{equation}
Taking $y=u-b$ and $g(y)=\phi_2(y+b)$ we can rewrite (\ref{phi 2 renewal equation}) as follows:
\begin{equation*}
    g(y)=\frac{\lambda r}{c_2}\int_0^{y}g(y-x)T_{\rho_2}f(x)dx+\frac{\lambda r}{c_2}h(y), \ \ \ \ \ \ \ \ \ \ \ \  y>0,
\end{equation*}
where
\begin{equation*}
    h(y)=h(u-b)=\int_{u-b}^u\phi_1(u-x)T_{\rho_2}f(x)dx+T_{\rho_2}\bar{F}(u),\ \ u>b.
\end{equation*}
Hence
\begin{eqnarray*}
  \phi_2(u)&=& g(y) \\
           &=& \frac{\lambda r}{c_2}\int_0^{y}g(y-x)T_{\rho_2}f(x)dx+\frac{\lambda r}{c_2}h(y) \\
           &=& \sum_{n=0}^{\infty}\left(\frac{\lambda r}{c_2}\right)^{n+1} (T_{\rho_2}f)^{n*}*h(y)\\
           &=& \sum_{n=0}^{\infty}\left(\frac{\lambda r}{c_2}\right)^{n+1} (T_{\rho_2}f)^{n*}*h(u-b)
\end{eqnarray*}
which completes the proof.
\vspace{3mm} \hfill $\Box$

The expression for $\phi_1(u)$ could be also derived in terms of the translation operator.

\begin{theo}
The function $\phi_1(u)$ can be expressed analytically in the following form:
\begin{equation}\label{renewal equation 1}
   \phi_1(u)=\phi_{\infty}(u)+\frac{\frac{\lambda r}{c_2}\left[\phi_{\infty}*T_{\rho_2}f(b)+T_{\rho_2}\bar{F}(b)\right]-\phi_{\infty}(b)}{\nu(b)-\frac{\lambda r}{c_2}\nu*T_{\rho_2}f(b)}\nu(u),
\end{equation}
where
\begin{equation}\label{phi infty}
    \phi_{\infty}(u):=\sum_{n=0}^{\infty}\left(\frac{\lambda r}{c_1}\right)^{n+1}(T_{\rho_1} f)^{n*}*T_{\rho_1}\bar{F}(u)
\end{equation}
and
\begin{equation}\label{change2}
\nu(x):=\sum_{n=0}^{\infty}\left(\frac{\lambda r}{c_1}\right)^{n}(T_{\rho_1} f)^{n*}*p(x)\end{equation}
with $p(x)=e^{\rho_1x}$.
\end{theo}

\textbf{Proof.} Note that the first equation in (\ref{phi'}) does not involve the barrier level $b$:
\begin{equation}\label{extended_phi_11'}
    \phi_1'(u)=\frac{\lambda+\delta}{c_1}\phi_1(u)-\frac{\lambda r}{c_1}\int_0^u\phi_1(u-x)f(x)dx-\frac{\lambda r}{c_1}\bar{F}(u).
\end{equation}
The information about the barrier $b$ is included in the boundary condition:
$$\phi_1(b)=\phi_2(b):=\lim_{u\rightarrow b^{+}}\phi_2(u).$$
Lin et al. \cite{LiPa2006} showed that the general solution of (\ref{extended_phi_11'}) is of the form
\begin{equation}\label{phi 1 jie}
   \phi_1(u)=\phi_{\infty}(u)+k\nu(u),
\end{equation}
where $\phi_{\infty}(u)$ is the joint Laplace transform of density of the ruin time and number of claims counted up to ruin time
for the classical risk process (\ref{CL}) without any barrier applied.
That is,
\begin{equation}\label{w02}
    \phi_{\infty}(u):=\sum_{n=1}^\infty r^n\int_0^\infty e^{-\delta t}w_{\infty}(u,n,t)dt
\end{equation}
for
\begin{equation}\label{winfty}
w_{\infty}(u,n,t):=P(N_\tau=n, \tau\in dt|U(0)=u)/dt.\end{equation}

In above equation (\ref{phi 1 jie}) the quantity $k$ is a constant which we can specify by implementing (\ref{phi 1 jie})
and (\ref{phi 2 renewal equation}):
\begin{equation}\label{k}
    k=\frac{\frac{\lambda r}{c_2}\left[\int_{0}^b\phi_{\infty}(b-x)T_{\rho_2}f(x)dx+T_{\rho_2}\bar{F}(b)\right]-\phi_{\infty}(b)}{\nu(b)-\frac{\lambda r}{c_2}\int_{0}^b\nu(b-x)T_{\rho_2}f(x)dx}.
\end{equation}

We express now the function $\phi_{\infty}$ in terms of a compound geometric distribution.
Indeed, since $\phi_{\infty}$ also satisfies equation (\ref{extended_phi_11'}), taking Laplace transforms of its both sides for sufficiently large $s$
gives:
   \begin{equation}\label{Laplace_phi_11}
    (c_1s-(\lambda+\delta)+\lambda r\hat{f}(s))\hat{\phi}_{\infty}(s)=c_1\phi_{\infty}(0)-\lambda r\hat{\bar{F}}(s),\ s\geq0.
\end{equation}
To determine the constant term $c_1\phi_{\infty}(0)$ in (\ref{Laplace_phi_11}), we substitute the solution $\rho_1$
of the Lundberg Fundamental Equation (\ref{Lundberg_fundamental_equation}) for $s$:
\begin{equation}\label{phi0}
    c_1\phi_{\infty}(0)=\lambda r\hat{\bar{F}}(\rho_1) =\lambda r T_{\rho_1}\hat{\bar{F}}(0).
\end{equation}
Consequently, the equation (\ref{Laplace_phi_11}) reduces to
\begin{equation*}
    [c_1(s-\rho_1)+\lambda r\hat{f}(s)-\lambda r\hat{f}(\rho_1)]\hat{\phi}_{\infty}(s)=\lambda r\hat{\bar{F}}(\rho_1)-\lambda r\hat{\bar{F}}(s).
   \end{equation*}
Dividing above equation by $s-\rho_1$ and simple rearranging along with implementation of the formula (\ref{Dickson's Operator}) produces:
\begin{equation*}
    c_1\hat{\phi}_{\infty}(s)=\lambda r\hat{\phi}_{\infty}(s)T_sT_{\rho_1}f(0)+\lambda rT_sT_{\rho_1}\bar{F}(0).
\end{equation*}
Inverting this Laplace transforms gives classical renewal equation:
\begin{equation}\label{renewal equation phi_1}
   \phi_{\infty}(u)=\frac{\lambda r}{c_1}\phi_{\infty}*T_{\rho_1}f(u)+\frac{\lambda r}{c_1}T_{\rho_1}\bar{F}(u)
\end{equation}
having the solution given as an Neumann infinite series (\ref{phi infty}).

To prove the last statement (\ref{change2}) note that the function $\nu(u)$ satisfies the following integro-differential equation:
\begin{equation}\label{jian1}
   c_1\nu^{'}(u)-(\lambda+\delta)\nu(u)+\lambda r\int_0^{u}\nu(u-x)f(x)dx=0,  \ \ \ u \geq0,
\end{equation}
with the initial condition $\nu(0)=1$.
To get the analytical expression of $\nu(u)$
we take the Laplace transforms of both sides of (\ref{jian1}) for sufficiently large $s$ ($s>\rho_1$). This yields:
\begin{equation*}
    c_1s\hat{\nu}(s)-c_1\nu(0)=(\lambda+\delta)\hat{\nu}(s)-\lambda r\hat{f}(s)\hat{\nu}(s).
\end{equation*}
Since $\nu(0)=1$,
\begin{equation}\label{lundberg nu}
    (s+\frac{\lambda r}{c_1}\hat{f}(s)-\frac{\lambda+\delta}{c_1})\hat{\nu}(s)=1.
\end{equation}
Recalling that $\rho_1$ is the root of (\ref{Lundberg_fundamental_equation}), we can rewrite (\ref{lundberg nu}) as
\begin{equation*}\label{change}
    (s-\rho_1+\frac{\lambda r}{c_1}[\hat{f}(s)-\hat{f}(\rho_1)])\hat{\nu}(s)=1,
\end{equation*}
which, by dividing by $s-\rho_1$ and implementing (\ref{Dickson's Operator}), produces:
\begin{equation}\label{change1}
   \hat{\nu}(s)=\frac{\lambda r}{c_1}\hat{\nu}(s)T_sT_{\rho_1}f(0)+\frac{1}{s-\rho_1}.
\end{equation}
Inverting the Laplace transforms in (\ref{change1}) leads to the equation (\ref{change2}).
Including all above identities in (\ref{phi 1 jie}) completes the proof.
\vspace{3mm} \hfill $\Box$

\section{The joint density of $(\tau, N_\tau)$ }
In this section we
give the joint density of the number of claims counted until ruin time and the ruin time using
the Lagrange's Expansion theorem.
We start with few facts
that will be useful in the proof of the main result.

Recall that by $w_{\infty}(u,n,t)$ we denote the joint density of $(\tau, N_\tau)$ for the classical risk process (\ref{CL})
(with infinite barrier $b=+\infty$); see (\ref{winfty}).
For $i=1,2$ we denote
\begin{eqnarray*}
g_{i}(x,0,t) &:=& \delta_{{x}/{c_i}}(t)e^{-\lambda x/c_i}, \\
g_{i}(x,n,t) &:=& x{t^{n-1}}e^{-\lambda t}\lambda^nf^{n*}(c_it-x)/{n!}.
\end{eqnarray*}

\begin{lem}\label{winftylem}
We have
\begin{equation*}\label{w infty(u,1,t)}
     w_{\infty}(u,1,t)=\lambda e^{-\lambda t}\bar{F}(u+c_1t).
 \end{equation*}
For $n=1,2,3,\ldots$ the following holds:
 \begin{eqnarray}\label{w infty(u,n,t)}
   w_{\infty}(u,n+1,t)&=&\frac{({\lambda t})^{n}}{n!}e^{-\lambda t}\int_0^{u+c_1t}f^{n*}(u+c_1t-x)\lambda\bar{F}(x)dx\nonumber\\
    &&-c_1\sum_{j=1}^{n}\int_0^t\frac{({\lambda s})^{j}}{j!}e^{-\lambda s}f^{j*}(u+c_1s)w_{\infty}(0,n+1-j,t-s))ds,
 \end{eqnarray}
where
\begin{equation}\label{w infty(0,n,t)}
    w_{\infty}(0,n,t)=\frac{\lambda}{c_1}\int_0^{c_1t} \bar{F}(x)g_{1}(x,n-1,t)dx,\ \ \ n=1,2,\ldots.
\end{equation}
\end{lem}
\textbf{Proof.}
Using Lagrange's Expansion Theorem presented in Section 2 with
$\alpha(z)=e^{-zx}$, $\beta(z)=-\frac{\lambda}{c_i}\hat{f}(s)$, $a=(\lambda+\delta)/c_i$ and $D=\{z||z-a|\leq a\}$ ($i=1,2$) and
the Lundberg fundamental equation (\ref{Lundberg_fundamental_equation}) we can conclude the following identity:
\begin{eqnarray*}
e^{-\rho_i x} &=& e^{-(\lambda+\delta)x/c_i}+\sum_{n=1}^\infty\frac{r^n}{n!}\frac{d^{n-1}}{ds^{n-1}}\left(-xe^{-sx}\left(-\frac{\lambda}{c_i}\hat{f}(s)\right)^n\right)\Big|_{s=(\lambda+\delta)/c_i}\nonumber\\
&=&e^{-(\lambda+\delta)x/c_i}+\sum_{n=1}^\infty\frac{r^n}{n!}\frac{d^{n-1}}{ds^{n-1}}\left((-1)^{n+1}\lambda^nx/{c_i^n} \int_0^\infty e^{-s(x+y)}f^{n*}(y)dy\right)\Big|_{s=(\lambda+\delta)/c_i}\nonumber\\
&=&e^{-(\lambda+\delta)x/c_i}+\sum_{n=1}^\infty\frac{\lambda^nr^n}{n!c_i^n}\int_0^\infty x(x+y)^{n-1}e^{-(\lambda+\delta)(x+y)/c_i}f^{n*}(y)dy.
\end{eqnarray*}
Substituting $t:=(x+y)/c_i$ and rearranging leads to:
\begin{eqnarray}\label{e rho 1}
e^{-\rho_i x} &=& e^{-(\lambda+\delta)x/c_i}+\sum_{n=1}^\infty r^n\frac{\lambda^n}{n!}\int_{x/c_i}^{\infty} x t^{n-1}e^{-\lambda t}e^{-\delta t}f^{n*}(c_it-x)dt\nonumber\\
&=&\sum_{n=0}^\infty r^n\int_{x/c_i}^{\infty} e^{-\delta t}g_{i}(x,n,t)dt.
\end{eqnarray}
Therefore,
\begin{eqnarray}\label{T rho 1 f}
  T_{\rho_i}f(x) &=& \int_x^\infty e^{-\rho_i(u-x)}f(u)du\nonumber \\
                 &=& \int_x^\infty\sum_{n=0}^\infty r^n\int_{(u-x)/c_i}^{\infty} e^{-\delta t}g_{i}(u-x,n,t)dt f(u)du \nonumber\\
                 &=& \sum_{n=0}^\infty r^n \int_0^\infty e^{-\delta t}\int_x^{c_it+x}f(u)g_{i}(u-x,n,t)du dt.
\end{eqnarray}
Since $\phi_{\infty}(u)$ defined in (\ref{w02})
is the joint Laplace transform under the classical compound Poisson risk model without a barrier we can use
Dickson \cite{Di2012} to complete the proof.
\vspace{3mm} \hfill $\Box$

\begin{lem}
The function $\nu(u)$ given in (\ref{change2}) equals
\begin{equation}\label{nu density}
    \nu(u)=\sum_{n=0}^\infty r^n \int_0^\infty e^{-\delta t}\varpi(u,n,t)dt,
\end{equation}
where
\begin{eqnarray*}
 \varpi(u,0,t)&:=&g_{1}(-u,0,t),\\
 \varpi(u,n,t)&:=&\sum_{m=1}^n\left(\frac{\lambda}{c_1}\right)^m\int_0^{c_1t}\int_0^u g_{c_1}(y,n-m,t) b_m(u-x,y+x)dx dy+g_{c_1}(-u,n,t),\ \ \ n\geq1\\
  b_n(u,y) &:=& \sum_{j=0}^{n-1}{n \choose j}\frac{(-1)^j}{\Gamma(n)}\int_0^u (u-x)^{n-1}f^{(n-j)*}(y+u-x)f^{j*}(x)dx.
\end{eqnarray*}
\end{lem}

\textbf{Proof.}
Our goal is to express $\nu(u)$ as the Laplace transform:
\begin{eqnarray}\label{change2+}
  \nu(u) &=& \int_0^\infty e^{-\rho_1t}\xi(u,t)dt.
\end{eqnarray}

We start from definition (\ref{change2}):
\begin{eqnarray}\label{change2++}
  \nu(u) &=& \sum_{n=0}^{\infty}\left(\frac{\lambda r}{c_1}\right)^{n}(T_{\rho_1} f)^{n*}*p(u)\nonumber \\
         &=& \sum_{n=1}^{\infty}\left(\frac{\lambda r}{c_1}\right)^{n}\int_0^u (T_{\rho_1} f)^{n*}(u-x)e^{\rho_1x}dx+e^{\rho_1u}.
\end{eqnarray}
Using Dickson and Willmot \cite{DiWi2005} we can obtain the following representation:
\begin{equation}\label{b n}
    (T_{\rho_i} f)^{n*}(u)=\int_0^\infty e^{-\rho_iy}b_n(u,y)dy
\end{equation}
for
\begin{equation*}
    b_n(u,y):=\sum_{j=0}^{n-1}{n \choose j}\frac{(-1)^j}{\Gamma(n)}\int_0^u (u-x)^{n-1}f^{(n-j)*}(y+u-x)f^{j*}(x)dx.
\end{equation*}
By (\ref{change2++})
\begin{eqnarray*}
  \nu(u) &=& \sum_{n=1}^{\infty}\left(\frac{\lambda r}{c_1}\right)^{n}\int_0^u \int_0^\infty e^{-\rho_1y}b_n(u-x,y)dye^{\rho_1x}dx+e^{\rho_1u} \\
   &=& \sum_{n=1}^{\infty}\left(\frac{\lambda r}{c_1}\right)^{n}\int_0^\infty e^{-\rho_1t} \int_0^u b_n(u-x,t+x)dx dt \\
   && +\sum_{n=1}^{\infty}\left(\frac{\lambda r}{c_1}\right)^{n}\int_{-u}^0 e^{-\rho_1t} \int_{-t}^u b_n(u-x,t+x)dx dt+\int_0^\infty e^{-\rho_1t}\delta_{-u}(t)dt.
\end{eqnarray*}
Comparing the coefficients of $e^{-\rho_1t}$ in (\ref{change2+}) gives:
\begin{equation}\label{xirow}
    \xi(u,t)=\sum_{n=1}^{\infty}\left(\frac{\lambda r}{c_1}\right)^{n}\int_0^u b_n(u-x,t+x)dx+\delta_{-u}(t).
\end{equation}
Using (\ref{e rho 1}) and (\ref{xirow}) in (\ref{change2+}) we end up with:
\begin{eqnarray*}
  \nu(u) &=& \int_0^\infty e^{-\rho_1y}\xi(u,y)dy+e^{\rho_1u} \\
   &=& \int_0^\infty \sum_{n=0}^\infty r^n\int_{y/c_1}^{\infty} e^{-\delta t}g_{c_1}(y,n,t)dt\xi(u,y)dy \\
   &=& \sum_{n=0}^\infty r^n \int_0^\infty e^{-\delta t}\int_0^{c_1t}g_{c_1}(y,n,t)\xi(u,y)dy dt\\
   &=& \sum_{n=1}^\infty r^n \int_0^\infty e^{-\delta t}\left(\sum_{m=1}^n(\frac{\lambda}{c_1})^m\int_0^{c_1t}\int_0^u g_{c_1}(y,n-m,t) b_m(u-x,y+x)dx dy+g_{c_1}(-u,n,t)\right)dt\\
   &&+\int_0^\infty e^{-\delta t}g_{c_1}(-u,0,t)dt
\end{eqnarray*}
which completes the proof.
\vspace{3mm} \hfill $\Box$

Using  above lemmas we will prove the main result of this paper.
\begin{theo}
For $0\leq u\leq b$ and $m>1$ the joint density of the number of claims until ruin $N_\tau$ and the time to ruin $\tau$ is given by
\begin{eqnarray}\label{w(u,n,t)}
     &&w_1(u,1,t)=\frac{\lambda}{c_2}e^{-\lambda t}\bar{F}(c_2t+b+\frac{c_2}{c_1}(u-b))\nonumber\\
     &&w_1(u,m,t)=e^{-\frac{\lambda b}{c_1}}\left[\sum_{n=1}^m\vartheta(u,m,n,t-\frac{b}{c_1})-\sum_{n=1}^{m-1}\int_0^{t-\frac{b}{c_1}}\varsigma(b,m-n,t-\frac{b}{c_1}-z)w_1(u,n,z)dz\right],\nonumber\\
\end{eqnarray}
where for $n\geq 1$
\begin{eqnarray*}
  &&\varsigma(b,0,t) := \varpi(b,0,t)= g_{c_1}(-b,0,t), \\
  &&\varsigma(b,n,t) := \varpi(b,n,t)-\sum_{m=0}^{n-1}\frac{\lambda }{c_2}\int_0^b\int_0^t\varpi(b-x,n-1-m,t-z)\int_x^{c_2z+x}f(y)g_{2}(y-x,m,z)dy dz dx,\\
  &&\gamma(b,1,t) :=\frac{\lambda}{c_2}\int_b^{c_2t+b}\bar{F}(y)g_{2}(y-b,0,t)dy,  \\
  &&\gamma(b,n,t) := \sum_{m=0}^{n-2}\frac{\lambda }{c_2}\int_0^b\int_0^tw_\infty(b-x,n-m-1,t-z)\int_x^{c_2z+x}f(y)g_{2}(y-x,m,z)dy dz dx,\\
  &&\vartheta(u,m,n,t):=\int_0^t\varsigma(b,m-n,t-z)w_\infty(u,n,z)+\left(\gamma(b,n,t-z)-w_\infty(b,n,t-z)\right)\varpi(u,m-n,z)dz.
\end{eqnarray*}
\end{theo}
\textbf{Proof.}
In order to get the joint density $w(u,n,t)$, we have to take inverse Laplace transform with respect to $\delta$ rather than $\rho_1$ and $\rho_2$.
To do this we must find firstly the relationship between transforms with respect to $\rho_1$, $\rho_2$ and $\delta$ by applying the Lagrange's Expansion theorem.
For convenience, we will denote:
\begin{equation}\label{chi}
    \chi(b):=\nu(b)-\frac{\lambda r}{c_2}\nu*T_{\rho_2}f(b).
\end{equation}
Then we can rewrite (\ref{renewal equation 1}) as follows:
\begin{equation}\label{456}
    \chi(b)\phi_1(u)=\chi(b)\phi_{\infty}(u)+\frac{\lambda r}{c_2}\left[\phi_{\infty}*T_{\rho_2}f(b)+T_{\rho_2}\bar{F}(b)\right]\nu(u)-\phi_{\infty}(b)\nu(u).
\end{equation}
Putting (\ref{T rho 1 f}) and (\ref{nu density}) into (\ref{chi}) we will derive:
\begin{eqnarray}\label{chi1}
  \chi(b) &=& \sum_{n=0}^\infty r^n \int_0^\infty e^{-\delta t}\varpi(b,n,t)dt-\frac{\lambda r}{c_2}\int_0^b\nu(b-x)T_{\rho_2}f(x)dx\nonumber \\
   &=& \sum_{n=0}^\infty r^n \int_0^\infty e^{-\delta t}\varpi(b,n,t)dt-\frac{\lambda r}{c_2}\sum_{n=0}^\infty r^n\sum_{m=0}^n\int_0^b\int_0^\infty e^{-\delta t}\varpi(b-x,n-m,t)dt\int_0^\infty e^{-\delta z}\nonumber \\
   && \int_x^{c_2z+x}f(y)g_{2}(y-x,m,z)dy dz dx\nonumber\\
   &=& \sum_{n=0}^\infty r^n \int_0^\infty e^{-\delta t}\varpi(b,n,t)dt-\sum_{n=1}^\infty r^n \int_0^\infty e^{-\delta t}\{\sum_{m=0}^{n-1}\frac{\lambda }{c_2}\int_0^b\int_0^t\varpi(b-x,n-1-m,t-z)\nonumber\\
   && \int_x^{c_2z+x}f(y)g_{2}(y-x,m,z)dy dz dx\}dt\nonumber\\
   &=& \sum_{n=1}^\infty r^n \int_0^\infty e^{-\delta t}\{\varpi(b,n,t)-\sum_{m=0}^{n-1}\frac{\lambda }{c_2}\int_0^b\int_0^t\varpi(b-x,n-1-m,t-z)\int_x^{c_2z+x}f(y)\nonumber\\
   &&g_{2}(y-x,m,z)dy dz dx\}dt+\int_0^\infty e^{-\delta t}\varpi(b,0,t)dt\nonumber\\
   &=& \sum_{n=0}^\infty r^n \int_0^\infty e^{-\delta t}\varsigma(b,n,t)dt.
\end{eqnarray}
Similarly, using Lemma \ref{winftylem}, we can check that:
\begin{equation}\label{part}
    \frac{\lambda r}{c_2}\left[\phi_{\infty}*T_{\rho_2}f(b)+T_{\rho_2}\bar{F}(b)\right]=\sum_{n=1}^\infty r^n \int_0^\infty e^{-\delta t}\gamma(b,n,t)dt.
\end{equation}
Using (\ref{nu density}), (\ref{chi1}) and (\ref{part}) in (\ref{456}) we obtain:
\begin{eqnarray*}
    &&\sum_{m=1}^\infty r^m \int_0^\infty e^{-\delta t}\sum_{n=1}^m\int_0^t\varsigma(b,m-n,t-z)\left(w_1(u,n,z)-w_\infty(u,n,z)\right)dz dt\\
    &=&\sum_{m=1}^\infty r^m \int_0^\infty e^{-\delta t}\sum_{n=1}^m\int_0^t\left(\gamma(b,n,t-z)-w_\infty(b,n,t-z)\right)\varpi(u,m-n,z)dz dt
\end{eqnarray*}
or equivalently that
\begin{eqnarray*}
 \sum_{n=1}^m\int_0^t\varsigma(b,m-n,t-z)w_1(u,n,z)dz &=&\sum_{n=1}^m\int_0^t\varsigma(b,m-n,t-z)w_\infty(u,n,z)\\
 &&+\left(\gamma(b,n,t-z)-w_\infty(b,n,t-z)\right)\varpi(u,m-n,z)dz.
\end{eqnarray*}

Now, if $m=1$ then
\begin{eqnarray*}
  \int_0^t\varsigma(b,0,t-z)w_1(u,1,z)dz &=&  \vartheta(u,1,1,t).
\end{eqnarray*}
In this case
\begin{eqnarray*}
  \int_0^t\delta_{-b/c_1}(t-z)e^{\frac{\lambda b}{c_1}}w_1(u,1,z)dz &=& e^{\frac{\lambda b}{c_1}}w_1(u,1,t+\frac{b}{c_1}) \\
                                                                    &=& \frac{\lambda}{c_2}e^{-\lambda t}\bar{F}(c_2t+b+\frac{c_2}{c_1}u)
\end{eqnarray*}
and
\begin{equation*}
    w_1(u,1,t)=\frac{\lambda}{c_2}e^{-\lambda (t+\frac{b}{c_1})}\bar{F}(c_2t+b+\frac{c_2}{c_1}(u-b)).
\end{equation*}
Similarly, if $m=2$ then
\begin{eqnarray*}
 \int_0^t\delta_{-b/c_1}(t-z)e^{\frac{\lambda b}{c_1}}w_1(u,2,z)dz &=& \sum_{n=1}^2\vartheta(u,2,n,t)-\int_0^t\varsigma(b,1,t-z)w_1(u,1,z)dz
\end{eqnarray*}
and
\begin{equation*}
     w_1(u,2,t)=e^{-\frac{\lambda b}{c_1}}\left[\sum_{n=1}^2\vartheta(u,2,n,t-\frac{b}{c_1})-\int_0^{t-\frac{b}{c_1}}\varsigma(b,1,t-\frac{b}{c_1}-z)w_1(u,1,z)dz\right].
\end{equation*}
Similarly we can prove the assertion for any $m>1$.
\vspace{3mm} \hfill $\Box$

\begin{theo}
For $u >b$ and $m>1$ the joint density of the number of claims until ruin $N_\tau$ and the time to ruin $\tau$ is given by
\begin{equation}\label{w(u,n,t)}
    w_2(u,m,t)=(\frac{\lambda}{c_2})^{m}\sum_{k=0}^{m-1}\sum_{n=0}^{m-k-1}\int_0^{u-b}\int_0^t\int_0^{c_2z}
   g_{2}(y,k,z) b_{m-k-n-1}(u-b-x,y)\varepsilon(x,n,t-z)dy dz dx,
\end{equation}
where
\begin{eqnarray*}
  \varepsilon(u,0,t) &:=& \int_{u+b}^{c_2t+u+b}\bar{F}(y)g_{2}(y-u-b,0,t)dy, \\
  \varepsilon(u,m,t) &:=& \sum_{n=1}^m\int_{u}^{u+b}\int_0^tw_1(u+b-x,n,t-z)\int_x^{c_2z+x}f(y)g_{2}(y-x,m-n,z)dy dz dx\nonumber\\
  &&+\int_{u+b}^{c_2t+u+b}\bar{F}(y)g_{2}(y-u-b,m,z)dy, \ \ n\geq1.
\end{eqnarray*}
\end{theo}

\textbf{Proof.}
To obtain an expression for $w_2(u,m,t)$ we first consider $h(x)$ defined in (\ref{h}).
Using (\ref{T rho 1 f}) we can derive:
\begin{eqnarray}\label{h x}
  h(u) &=& \int_{u}^{u+b}\sum_{m=1}^\infty r^m \int_0^\infty e^{-\delta t}w_1(u+b-x,m,t)dt\nonumber\\&&\quad\times\sum_{n=0}^\infty r^n \int_0^\infty e^{-\delta z}\int_x^{c_2z+x}f(y)g_{2}(y-x,n,z)dy dz dx\nonumber\\
  &&+ \sum_{n=0}^\infty r^n \int_0^\infty e^{-\delta t}\int_{u+b}^{c_2t+u+b}\bar{F}(y)g_{2}(y-u-b,n,t)dy dt\nonumber\\
  &=& \sum_{n=1}^\infty r^n \int_0^\infty e^{-\delta t}\Biggl[\sum_{m=1}^n\int_{u}^{u+b}\int_0^tw_1(u+b-x,m,t-z)\nonumber\\&&\quad\times\int_x^{c_2z+x}f(y)g_{2}(y-x,n-m,z)dy dz dx\nonumber\Biggr.\\
  &&+\Biggl.\int_{u+b}^{c_2t+u+b}\bar{F}(y)g_{2}(y-u-b,n,z)dy\Biggr] dt+\int_0^\infty e^{-\delta t}\int_{u+b}^{c_2t+u+b}\bar{F}(y)g_{2}(y-u-b,0,t)dy dt\nonumber\\
  &=& \sum_{n=0}^\infty r^n \int_0^\infty e^{-\delta t}\varepsilon(u,n,t)dt.
\end{eqnarray}
Moreover, substituting (\ref{b n}), (\ref{h x}) and (\ref{e rho 1}) into (\ref{renewal equation 2}) gives:
\begin{eqnarray}\label{phi 2 u}
  \phi_2(u) &=& \sum_{m=0}^{\infty}\left(\frac{\lambda r}{c_2}\right)^{m+1}\int_0^{u-b}(T_{\rho_2}f)^{m*}(u-b-x)h(x)dx\nonumber \\
   &=& \sum_{m=0}^{\infty}\left(\frac{\lambda r}{c_2}\right)^{m+1}\int_0^{u-b}\int_0^\infty e^{-\rho_2y}b_m(u-b-x,y)dy\sum_{n=0}^\infty r^n \int_0^\infty e^{-\delta t}\varepsilon(x,n,t)dtdx\nonumber \\
   &=& \sum_{m=0}^{\infty}\left(\frac{\lambda r}{c_2}\right)^{m+1}\int_0^{u-b}\int_0^\infty \sum_{k=0}^\infty r^k\int_{y/c_2}^{\infty} e^{-\delta z}g_{2}(y,k,z)dz b_m(u-b-x,y)dy\nonumber\\
   &&\times\sum_{n=0}^\infty r^n \int_0^\infty e^{-\delta t}\varepsilon(x,n,t)dtdx\nonumber\\
   &=&\sum_{m=1}^\infty r^m \int_0^\infty e^{-\delta t}\Biggl\{\left(\frac{\lambda}{c_2}\right)^{m}\sum_{k=0}^{m-1}\sum_{n=0}^{m-k-1}\int_0^{u-b}\int_0^t\int_0^{c_2z}
   g_{2}(y,k,z) b_{m-k-n-1}(u-b-x,y)\Biggr.\nonumber\\
   &&\Biggl.\times\varepsilon(x,n,t-z)dy dz dx\Biggr\}dt.
\end{eqnarray}
Comparing equations (\ref{phi 2 u}) and (\ref{define_phi}) completes the proof.
\vspace{3mm} \hfill $\Box$

\noindent{\bf{Acknowledgements}}\\
This research is  support by the National Natural Science Foundation of China (Grant No. 11271164)
and by the FP7 Grant PIRSES-GA-2012-318984.
Zbigniew Palmowski is supported by the National Science Centre under the grant 2015/17/B/ST1/01102
(2016-2019).

\bibliographystyle{plain}

\begin{thebibliography}{99}

\bibitem{BaDrLa2007}
A.~Badescu, S.~Drekic, and D.~Landriault.
\newblock Analysis of a threshold dividend strategy for a {M}{A}{P} risk model.
\newblock {\em Scandinavian Actuarial Journal}, 2007(4):227--247, 2007.

\bibitem{De1957}
B.~De~Finetti.
\newblock Su un{'}impostazione alternativa della teoria collettiva del rischio.
\newblock In {\em Transactions of the XVth international congress of
  Actuaries}, pages 433--443, 1957.

\bibitem{Di2012}
D.~C.~M. Dickson.
\newblock The joint distribution of the time to ruin and the number of claims
  until ruin in the classical risk model.
\newblock {\em Insurance: Mathematics and Economics}, 50(3):334--337, 2012.

\bibitem{DiHi2001}
D.~C.~M. Dickson and C.~Hipp.
\newblock On the time to ruin for {E}rlang (2) risk processes.
\newblock {\em Insurance: Mathematics and Economics}, 29(3):333--344, 2001.

\bibitem{DiWi2005}
D.~C.~M. Dickson and G.~E. Willmot.
\newblock The density of the time to ruin in the classical {P}oisson risk
  model.
\newblock {\em Astin Bulletin}, 35(1):45--60, 2005.

\bibitem{FrPiPo2012}
E~Frostig, S.~M. Pitts, and K.~Politis.
\newblock The time to ruin and the number of claims until ruin for phase-type
  claims.
\newblock {\em Insurance: Mathematics and Economics}, 51(1):19--25, 2012.

\bibitem{GaYi2008}
H.~Gao and C.~Yin.
\newblock The perturbed {S}parre {A}ndersen model with a threshold dividend
  strategy.
\newblock {\em Journal of Computational and Applied Mathematics},
  220(1-2):394--408, 2008.

\bibitem{GeSh2005}
H.~U. Gerber and E.~S.~W. Shiu.
\newblock The time value of ruin in a {S}parre {A}ndersen model.
\newblock {\em North American Actuarial Journal}, 9(2):49--69, 2005.

\bibitem{KyLo2010}
A.~E. Kyprianou and R.~L. Loeffen.
\newblock Refracted {L}{\'e}vy processes.
\newblock {\em Annales De L' Institut Henri Poincar{\'e} Probabilit{\'e}s Et
  Statistiques}, 46(1):24--44, 2010.

\bibitem{La1770}
J.~L. Lagrange.
\newblock {\em Nouvelle m{\'e}thode pour r{\'e}soudre les {\'e}quations
  litt{\'e}rales par le moyen des s{\'e}ries}.
\newblock Chez Haude et Spener, Libraires de la Cour \& de l'Acad¨¦mie royale,
  1770.

\bibitem{LaShWi2011}
D.~Landriault, T.~Shi, and G.~E. Willmot.
\newblock Joint densities involving the time to ruin in the {S}parre {A}ndersen
  risk model under exponential assumptions.
\newblock {\em Insurance: Mathematics and Economics}, 49(3):371--379, 2011.

\bibitem{LiGa2004a}
S.~Li and J.~Garrido.
\newblock On a class of renewal risk models with a constant dividend barrier.
\newblock {\em Insurance: Mathematics and Economics}, 35(3):691--701, 2004.

\bibitem{LiGa2004b}
S.~Li and J.~Garrido.
\newblock On ruin for the {E}rlang (n) risk process.
\newblock {\em Insurance: Mathematics and Economics}, 34(3):391--408, 2004.


\bibitem{LiLu2009}
S.~Li and Y.~Lu.
\newblock The distribution of total dividend payments in a {S}parre {A}ndersen
  model.
\newblock {\em Statistics and Probability Letters}, 79(9):1246--1251, 2009.

\bibitem{LiLu2016}
S.~Li and Y.~Lu.
\newblock On the time and the number of claims when the surplus drops below a certain level.
\newblock  {\em Scandinavian Actuarial Journal}, 5: 420--445, 2016.

\bibitem{LiLu2}
Y.~Lu and S.~Li.
\newblock The Markovian regime-switching risk model with a threshold dividend strategy.
\newblock {\em Insurance: Mathematics and Economics}, 44(2):296--303, 2009.

\bibitem{LiPa2006}
X.~S. Lin and K.~P. Pavlova.
\newblock The compound {P}oisson risk model with a threshold dividend strategy.
\newblock {\em Insurance: Mathematics and Economics}, 38(1):57--80, 2006.

\bibitem{LiWiDr2003}
X.~S. Lin, G.~E. Willmot, and S.~Drekic.
\newblock The classical risk model with a constant dividend barrier: analysis
  of the {G}erber-{S}hiu discounted penalty function.
\newblock {\em Insurance: Mathematics and Economics}, 33(3):551--566, 2003.

\bibitem{WiDi2003}
G.~E. Willmot and D.~C.~M. Dickson.
\newblock The {G}erber-{S}hiu discounted penalty function in the stationary
  renewal risk model.
\newblock {\em Insurance: Mathematics and Economics}, 32(3):403--411, 2003.

\bibitem{ZhZh2013}
C.~Zhao and C.~Zhang.
\newblock Joint density of the number of claims until ruin and the time to ruin
  in the delayed renewal risk model with {E}rlang (n) claims.
\newblock {\em Journal of Computational and Applied Mathematics}, 244:102--114,
  2013.

\bibitem{Zh2005}
X.~Zhou.
\newblock On a classical risk model with a constant dividend barrier.
\newblock {\em North American Actuarial Journal}, 9(4):95--108, 2005.

\bibitem{ZhYa2008}
J.~Zhu and H.~Yang.
\newblock Ruin theory for a {M}arkov regime-switching model under a threshold
  dividend strategy.
\newblock {\em Insurance: Mathematics and Economics}, 42(1):311--318, 2008.

\end{thebibliography}

\end{document}